\documentclass[a4paper,11pt]{article}
\usepackage{amsfonts}
\usepackage{amsmath}
\usepackage{amssymb}
\usepackage[mathscr]{eucal}
\usepackage{latexsym}
\usepackage{theorem}

\input{tcilatex}

\begin{document}

\begin{center}
\textbf{INVERSE NODAL PROBLEM FOR P-LAPLACIAN ENERGY-DEPENDENT
STURM-LIOUVILLE EQUATION}

Hikmet KEMALOGLU

Department of Mathematics, University of Firat, 23119, Elazig, TURKEY

hkoyunbakan@gmail.com
\end{center}

\begin{quote}
\textbf{Abstract:}{\footnotesize \ In this study, inverse nodal problem is
solved for }$p-${\footnotesize Laplacian Schr\"{o}dinger equation with
energy-dependent potential with the Drichlet conditions. Asymptotic
estimates of eigenvalues, nodal points and nodal lengths are given by using
Pr\"{u}fer substitution. Especially, an explicit formula for potential
function is given by using nodal lengths. Results are more general than
classical }$p-${\footnotesize \ Laplacian Sturm Liouville problem. For the
proofs, it is used the methods given in the references \cite{lav3}, \cite%
{Wang}.}

$\smallskip \linebreak $


\end{quote}

\QTP{Body Math}
$\smallskip $\linebreak \textbf{\large 1. Introduction}

Consider the following eigenvalue problem for $p$-Laplacian%
\begin{equation}
-\left( u^{\prime (p-1)}\right) ^{\prime }=(p-1)(\lambda -q(x))u^{(p-1)} 
\tag{1.1}
\end{equation}%
with the conditions%
\begin{equation}
u(0)=u(1)=0  \tag{1.2}
\end{equation}%
where $p\in \lbrack 1,\infty )$ is a constant, $\lambda $ is the eigenvalue
and $q\in L^{2}(0,1).$ Equation (1.1) is also known as one-dimensional $p-$%
Laplacian eigenvalue equation. Note that when $p=2,$ equation (1.1) becomes
Sturm-Liouville equation as%
\begin{equation*}
-u^{\prime \prime }+q(x)u=\lambda u
\end{equation*}%
and inverse problem was solved for (1.1),(1.2) in the references \cite%
{binding}, \cite{bind2},{\footnotesize \ }\cite{brovn}, \cite{dra}, \cite%
{dra2} \cite{lav3}, \cite{valter}, \cite{Wang}.

The determination of the form of a differential operator from spectral data
associated with it has enjoyed close attention from a number of authors in
last years. One of the these operators is Sturm - Liouville operator. In the
typical formulation of the inverse Sturm-Liouville problem one seeks to
recover both $q(x)$ and constants by giving the eigenvalues with another
piece of spectral data. These data can take several forms, leading to many
versions of the problem. Especially, the recent interest is a study by Hald
and McLaughlin \cite{hald}, \cite{mc} wherein the given spectral information
consists of a set of nodal points of eigenfunctions for the Sturm-Liouville
problems. These results were extended to the case of problems with
eigenparameter dependent boundary conditions by Browne and Sleeman \cite%
{browne}. On the other hand , Law et al \cite{lav}, Law and Yang \cite{lav2}%
, solved the inverse nodal problem of determining the smoothness of the
potential function $q(x)$ of the Sturm-Liouville problem by using nodal
data. In the past few years, the inverse nodal problem of Sturm-Liouville
problem has been investigated by several authors \cite{but}, \cite{hik}, 
\cite{yang}.

When $q=0,$ consider the problem

\begin{equation*}
-\left( u^{\prime (p-1)}\right) ^{\prime }=(p-1)\lambda u^{(p-1)}
\end{equation*}

\begin{equation*}
u(0)=u(1)=0.
\end{equation*}%
The eigenvalues of this problem were given as \cite{lav3} 
\begin{equation*}
\lambda _{n}=\left( n\pi _{p}\right) ^{p},n=1,2,3,...,
\end{equation*}%
where%
\begin{equation*}
\pi _{p}=2\dint\limits_{0}^{1}\frac{dt}{\left( 1-t^{p}\right) ^{1/p}}=\frac{%
2\pi }{p\sin \left( \dfrac{\pi }{p}\right) }
\end{equation*}%
and associated eigenfunction is denoted by $S_{p}(x).$ $S_{p}(x)$ and $%
S_{p}^{\prime }(x)$ are periodic functions satisfying the identity%
\begin{equation*}
\left[ S_{p}(x)\right] ^{p}+\left[ S_{p}^{\prime }(x)\right] ^{p}=1
\end{equation*}%
for arbitrary $x\in 
\mathbb{R}
.$ These functions known as generalized sine and cos functions and for $p=2$
become $\sin x$ and $\cos x$ \cite{lin}$.$

Now, we must present some further properties of $S_{p}(x)$ for deriving a
more detailed asymptotic formulas. These formulas are crucial in the
solution of our problem.

\textbf{Lemma 1.1. }\cite{lav3}

\textbf{a)} For $S_{p}^{\prime }\neq 0,$%
\begin{equation*}
(S_{p}^{\prime })^{\prime }=-\left\vert \frac{S_{p}}{S_{p}^{\prime }}%
\right\vert ^{p-2}.S_{p}
\end{equation*}

\textbf{b)} 
\begin{equation*}
\left( S_{p}S_{p}^{\prime (p-1)}\right) ^{\prime }=\left\vert S_{p}^{\prime
}\right\vert ^{p}-(p-1)S_{p}^{p}=1-p\left\vert S_{p}\right\vert
^{p}=(1-p)+p\left\vert S_{p}^{\prime }\right\vert ^{p}.
\end{equation*}

According to the Sturm-Liouville Theory, the zero set $x_{n}=\{x_{j}^{(n)}%
\}_{j=1}^{n}$ of the eigenfunction $u_{n}(x)$ corresponding to $\lambda _{n}$
is called the nodal set and $l_{j}^{n}=x_{j+1}^{n}-x_{j}^{n}$ is defined as
the nodal length of $u_{n}(x).$ Using the nodal datas, some uniqueness
results, reconstruction and stability of potential functions have been
solved by many authors \cite{lav3}, \cite{Wang}.

Consider $p-$Laplacian eigenvalue problem%
\begin{equation}
-\left( u^{\prime (p-1)}\right) ^{\prime }=(p-1)(\lambda ^{2}-q(x)-2\lambda
r(x))u^{(p-1)}  \tag{1.3}
\end{equation}%
with the Dirichlet conditions%
\begin{equation}
u(0)=u(1)=0  \tag{1.4}
\end{equation}%
and Neumann boundary conditions%
\begin{equation}
u^{\prime }(0)=u^{\prime }(1)=0  \tag{1.5}
\end{equation}%
where $p>1,$ $\lambda $ is eigenvalue and $q(x)\in L^{2}\left( 0,1\right)
,r(x)\in W_{2}^{1}\left( 0,1\right) .$

For $p=2$, equation (2.2) becomes%
\begin{equation}
-y^{\prime \prime }+[q(x)+2\lambda r(x)]y=\lambda ^{2}y.  \tag{1.6}
\end{equation}

This equation is known as diffusion equation or quadratic of differential
pencil. The eigenvalue equation (1.6) is of important both classical and
quantum mechanics. For example, such problems arise in solving the
Klein--Gordon equations, which describe the motion of massless particles
such as photons. Sturm--Liouville energy-dependent equations are also used
for modelling vibrations of mechanical systems in viscous media (see \cite%
{jav}). We note that in this type problem the spectral parameter $\lambda $
is related to the energy of the system, and this motivates the terminology
\textquotedblleft energy-dependent\textquotedblright\ used for the spectral
problem of the form (1.6).Inverse problem of quadratic pencil have been
solved by many authors in the references \cite{kasim}, \cite{hus}, \cite{hry}%
, \cite{hik2}, \cite{hik}, \cite{nabi}, \cite{wang}, \cite{yang}, \cite%
{yang2}.

As in $p-$Laplacian Sturm-Liouville problem, for $p(x)=r(x)=0$ eigenvalues
of the problem (1.3),(1.4) are given 
\begin{equation*}
\lambda _{n}=\left( n\pi _{p}\right) ^{p}
\end{equation*}%
and associated eigenfunction is denoted by $S_{p}(x,\lambda ).$ We will note
it by $S_{p}(x)$ briefly.

This paper is organized as follows: In section 2, we give asymptotic formula
for eigenvalues, nodal points and nodal lengths. In section 3, we give a
reconstruction formula for differential pencil by using nodal parameters.

\bigskip

{\large 2}\textbf{\large . Asypmtotic Estimates of Nodal Parameters}

In this section, we study the properties of eigenvalues of the $p-$
Laplacian operator (1.3) with the Drichlet conditions. For this, we must
introduce the Pr\"{u}fer substitution. One can easily obtain similiar
results for Neumann problem.

We define a modified Pr\"{u}fer substution%
\begin{eqnarray}
u(x) &=&c(x)S_{p}\left( \lambda ^{2/p}\theta (x)\right)  \TCItag{2.1} \\
u^{\prime }(x) &=&\lambda ^{2/p}c(x)S_{p}^{\prime }\left( \lambda
^{2/p}\theta (x)\right)  \notag
\end{eqnarray}%
or%
\begin{equation}
\frac{u^{\prime }(x)}{u(x)}=\lambda ^{2/p}\frac{S_{p}^{\prime }\left(
\lambda ^{2/p}\theta (x)\right) }{S_{p}\left( \lambda ^{2/p}\theta
(x)\right) }.  \tag{2.2}
\end{equation}%
Differentiating the (2.2) with respect to $x$ and applying Lemma 1.1. one
can obtain that%
\begin{equation}
\theta ^{\prime }(x)=1-\frac{q}{\lambda ^{2}}S_{p}^{p}-\frac{2}{\lambda }%
rS_{p}^{p}  \tag{2.3}
\end{equation}

\textbf{Theorem 2.1.} The eigenvalues $\lambda _{n}$ of the Dirichlet
problem (1.3),(1.4) are%
\begin{equation}
\lambda _{n}^{2/p}=n\pi _{p}+\frac{1}{p(n\pi _{p})^{p-1}}\dint%
\limits_{0}^{1}q(t)dt+\frac{2}{p(n\pi _{p})^{\frac{p-2}{p}}}%
\dint\limits_{0}^{1}r(t)dt+O\left( \frac{1}{n^{\frac{p+2}{p}}}\right) 
\tag{2.4}
\end{equation}

\textbf{Proof :} For the problem (1.3),(1.4), Let $\lambda =\lambda
_{n},\theta (0)=0$ and $\theta (1)=\dfrac{n\pi _{p}}{\lambda _{n}^{2/p}}.$
Firstly, we should integrate both sides of (2.3) on $[0,1]$%
\begin{equation*}
\dfrac{n\pi _{p}}{\lambda _{n}^{2/p}}=1-\frac{1}{\lambda _{n}^{2}}%
\dint\limits_{0}^{1}q(t)S_{p}^{p}(t)dt-\frac{2}{\lambda _{n}}%
\dint\limits_{0}^{1}r(t)S_{p}^{p}(t)dt
\end{equation*}%
using the identity%
\begin{equation*}
\frac{d}{dt}\left[ S_{p}\left( \lambda _{n}^{2/p}\theta (t)\right)
S_{p}^{\prime }\left( \lambda _{n}^{2/p}\theta (t)\right) ^{p-1}\right]
=\left( 1-p\left\vert S_{p}\left( \lambda _{n}^{2/p}\theta (t)\right)
\right\vert ^{p}\right) \lambda _{n}^{2/p}\theta ^{\prime }(t)
\end{equation*}%
and Lemma 1.1. (b), we get%
\begin{eqnarray}
\dfrac{n\pi _{p}}{\lambda _{n}^{2/p}} &=&1-\frac{1}{\lambda _{n}^{2}p}%
\dint\limits_{0}^{1}q(t)dt-\frac{2}{\lambda _{n}p}\dint\limits_{0}^{1}r(t)dt
\notag \\
&&+\frac{1}{\lambda _{n}^{2}p}\dint\limits_{0}^{1}\frac{q(t)}{\lambda
_{n}^{2/p}\theta ^{\prime }(t)}\frac{d}{dt}\left[ S_{p}\left( \lambda
_{n}^{2/p}\theta (t)\right) S_{p}^{\prime }\left( \lambda _{n}^{2/p}\theta
(t)\right) ^{p-1}\right] dt  \notag \\
&&+\frac{2}{\lambda _{n}p}\dint\limits_{0}^{1}\frac{r(t)}{\lambda
_{n}^{2/p}\theta ^{\prime }(t)}\frac{d}{dt}\left[ S_{p}\left( \lambda
_{n}^{2/p}\theta (t)\right) S_{p}^{\prime }\left( \lambda _{n}^{2/p}\theta
(t)\right) ^{p-1}\right] dt.  \TCItag{2.5}
\end{eqnarray}

Then, using integration by parts, we have%
\begin{eqnarray*}
\dint\limits_{0}^{1}\frac{q(t)}{\lambda _{n}^{2/p}\theta ^{\prime }(t)}\frac{%
d}{dt}\left[ S_{p}S_{p}^{\prime }{}^{p-1}\right] dt &=&-\lambda
_{n}^{-2/p}\dint\limits_{0}^{1}G\left( \lambda _{n}^{2/p}\theta (t)\right) 
\frac{d}{dt}\left( \frac{q(t)}{\theta ^{\prime }(t)}\right) dt \\
&=&O\left( \frac{1}{\lambda _{n}^{2/p}}\right) ,
\end{eqnarray*}%
where%
\begin{equation*}
G\left( \lambda _{n}^{2/p}\theta (x)\right) =S_{p}\left( \lambda
_{n}^{2/p}\theta (x)\right) S_{p}^{\prime }\left( \lambda _{n}^{2/p}\theta
(x)\right) ^{p-1}
\end{equation*}%
and when $x=0,1$%
\begin{equation*}
G\left( \lambda _{n}^{2/p}\theta (x)\right) =0.
\end{equation*}%
Similarly, one can show that%
\begin{equation*}
\dint\limits_{0}^{1}\frac{r(t)}{\lambda _{n}^{2/p}\theta ^{\prime }(t)}\frac{%
d}{dt}\left[ S_{p}S_{p}^{\prime }{}^{p-1}\right] dt=O\left( \frac{1}{\lambda
_{n}^{2/p}}\right) .
\end{equation*}%
Inserting these values in (2.5) and after some straightforward computations,
we obtain (2.4).

\medskip

\textbf{Theorem 2.2. }For the problem (1.3), (1.4), The nodal points
expansion satisfies

\begin{eqnarray*}
x_{j}^{n} &=&\frac{j}{n}+\frac{j}{pn^{p+1}\left( \pi _{p}\right) ^{p}}%
\dint\limits_{0}^{1}q(t)dt+\frac{2j}{pn\frac{2p-2}{p}\left( \pi _{p}\right) 
\frac{2p-2}{p}}\dint\limits_{0}^{1}r(t)dt+\frac{2}{\left( n\pi _{p}\right)
^{^{\frac{p}{2}}}}\dint\limits_{0}^{x_{j}^{n}}r(x)S_{p}^{p}dx \\
&&+\frac{1}{\left( n\pi _{p}\right) ^{^{p}}}\dint%
\limits_{0}^{x_{j}^{n}}q(x)S_{p}^{p}dx+O\left( \frac{j}{n^{\frac{3p+2}{p}}}%
\right) .
\end{eqnarray*}

\bigskip

\textbf{Proof:} Let $\lambda =\lambda _{n}$ and integrating (2.3) from $0$
to $x_{j}^{n},$ we have

\begin{equation*}
\dfrac{j.\pi _{p}}{\lambda _{n}^{2/p}}=x_{j}^{n}-\dint\limits_{0}^{x_{j}^{n}}%
\frac{2r(x)}{\lambda _{n}}S_{p}^{p}dx-\dint\limits_{0}^{x_{j}^{n}}\frac{q(x)%
}{\lambda _{n}^{2}}S_{p}^{p}dx.
\end{equation*}%
By using the estimates of eigenvalues as

\begin{equation*}
\dfrac{1}{\lambda _{n}^{2/p}}=\frac{1}{n\pi _{p}}+\frac{1}{p\left( n\pi
_{p}\right) ^{p+1}}\dint\limits_{0}^{1}q(t)dt+\frac{2}{p\left( n\pi
_{p}\right) \frac{3p-2}{p}}\dint\limits_{0}^{1}r(t)dt+O\left( \frac{1}{n^{%
\frac{3p+2}{p}}}\right) ,
\end{equation*}%
we obtain

\begin{eqnarray*}
x_{j}^{n} &=&\frac{j}{n}+\frac{j}{pn^{p+1}\left( \pi _{p}\right) ^{p}}%
\dint\limits_{0}^{1}q(t)dt+\frac{2j}{pn\frac{2p-2}{p}\left( \pi _{p}\right) 
\frac{2p-2}{p}}\dint\limits_{0}^{1}r(t)dt+\frac{2}{\left( n\pi _{p}\right)
^{^{\frac{p}{2}}}}\dint\limits_{0}^{x_{j}^{n}}r(x)S_{p}^{p}dx \\
&&+\frac{1}{\left( n\pi _{p}\right) ^{^{p}}}\dint%
\limits_{0}^{x_{j}^{n}}q(x)S_{p}^{p}dx+O\left( \frac{j}{n^{\frac{3p+2}{p}}}%
\right) .
\end{eqnarray*}

\textbf{Theorem 2.3.} As, $n\rightarrow \infty ,$%
\begin{equation}
l_{j}^{n}=\dfrac{\pi _{p}}{\lambda _{n}^{2/p}}+\frac{2}{p\lambda _{n}}%
\dint\limits_{0}^{1}r(t)dt+\frac{1}{p\lambda _{n}^{2}}\dint%
\limits_{0}^{1}q(t)dt+O\left( \frac{1}{\lambda _{n}^{\frac{4+p}{p}}}\right) 
\tag{2.6}
\end{equation}

\textbf{Proof :} For large $n\in 
\mathbb{N}
,$ integrating (2.3) on $[x_{j}^{n},x_{j+1}^{n}]$ and then 
\begin{equation*}
\dfrac{\pi _{p}}{\lambda _{n}^{2/p}}=l_{j}^{n}-\frac{2}{\lambda }%
\dint\limits_{x_{j}^{n}}^{x_{j+1}^{n}}r(t)S_{p}^{p}dt-\frac{1}{\lambda ^{2}}%
\dint\limits_{x_{j}^{n}}^{x_{j+1}^{n}}q(t)S_{p}^{p}dt
\end{equation*}%
or%
\begin{eqnarray}
\dfrac{\pi _{p}}{\lambda _{n}^{2/p}} &=&l_{j}^{n}-\frac{2}{p\lambda _{n}}%
\dint\limits_{x_{j}^{n}}^{x_{j+1}^{n}}r(t)dt-\frac{1}{p\lambda _{n}^{2}}%
\dint\limits_{x_{j}^{n}}^{x_{j+1}^{n}}q(t)dt+\frac{2}{\lambda _{n}p}%
\dint\limits_{x_{j}^{n}}^{x_{j+1}^{n}}\frac{1}{\lambda _{n}^{2/p}\theta
^{\prime }(t)}\frac{d}{dt}\left[ S_{p}S_{p}^{\prime }{}^{p-1}\right] r(t)dt+
\notag \\
&&\frac{1}{\lambda _{n}^{2}p}\dint\limits_{x_{j}^{n}}^{x_{j+1}^{n}}\frac{1}{%
\lambda _{n}^{2/p}\theta ^{\prime }(t)}\frac{d}{dt}\left[ S_{p}S_{p}^{\prime
}{}^{p-1}\right] q(t)dt.  \TCItag{2.7}
\end{eqnarray}%
By Lemma 1.1. and similar process of Theorem 2.1, we obtain that%
\begin{eqnarray*}
\dint\limits_{x_{j}^{n}}^{x_{j+1}^{n}}\frac{r(t)}{\lambda _{n}^{2/p}\theta
^{\prime }(t)}\frac{d}{dt}\left[ S_{p}S_{p}^{\prime }{}^{p-1}\right] dt
&=&-\dint\limits_{j\pi _{p}}^{(j+1)\pi _{p}}\left( \frac{q(t)}{\lambda
_{n}^{2/p}\theta ^{\prime }(t)}\right) ^{\prime }G(\tau )\frac{d\tau }{%
\lambda _{n}^{2/p}\theta ^{\prime }(t)} \\
&=&O\left( \frac{1}{\lambda _{n}^{4/p}}\right) ,
\end{eqnarray*}%
where $G(\tau )=S_{p}(\tau )S_{p}^{\prime }(\tau ){}^{(p-1)}$ and $\tau
=\lambda _{n}^{2/p}\theta (x).$ Similarly one can show that%
\begin{equation*}
\dint\limits_{x_{j}^{n}}^{x_{j+1}^{n}}\frac{q(t)}{\lambda _{n}^{2/p}\theta
^{\prime }(t)}\frac{d}{dt}\left[ S_{p}S_{p}^{\prime }{}^{p-1}\right]
dt=O\left( \frac{1}{\lambda _{n}^{4/p}}\right) .
\end{equation*}%
Inserting this value in (2.7), we obtain%
\begin{equation*}
l_{j}^{n}=\dfrac{\pi _{p}}{\lambda _{n}^{2/p}}+\frac{2}{p\lambda _{n}}%
\dint\limits_{x_{j}^{n}}^{x_{j+1}^{n}}r(t)dt+\frac{1}{p\lambda _{n}^{2}}%
\dint\limits_{x_{j}^{n}}^{x_{j+1}^{n}}q(t)dt+O\left( \frac{1}{\lambda _{n}^{%
\frac{4+p}{p}}}\right)
\end{equation*}

and by Theorem 2.1.%
\begin{equation*}
l_{j}^{n}=\dfrac{1}{n}+\frac{2}{p(n\pi _{p})^{p/2}}\dint%
\limits_{x_{j}^{n}}^{x_{j+1}^{n}}r(t)dt+\frac{1}{p(n\pi _{p})^{p}}%
\dint\limits_{x_{j}^{n}}^{x_{j+1}^{n}}q(t)dt+O\left( \frac{1}{n^{\frac{4+p}{p%
}}}\right) .
\end{equation*}

\medskip

{\large 3}\textbf{\large . Reconstruction of Potential Function in
Differential Pencil}

In this section, we give an explicit formula for potential function. The
method used in the proof of the theorem is similiar to classical
Sturm-Liouville problem \cite{lav3}, \cite{Wang}.

\textbf{Theorem 3.1.} Let $q\in L^{2}(0,1)$ and $r\in W_{2}^{1}(0,1).$ Then,%
\begin{equation*}
q(x)=\lim_{n\rightarrow \infty }p\lambda _{n}^{2}\left( \frac{\lambda
_{n}^{2/p}l_{j}^{n}}{\pi _{p}}-\frac{2r(x)}{p\lambda _{n}}-1\right)
\end{equation*}%
for $j=j_{n}(x)=\max \{j:x_{j}^{n}\leq x\}.$

\textbf{Proof :} Applying the mean value theorem for integral to (2.6), with
fixed $n,$ there exists $z\in (x_{j}^{n},x_{j+1}^{n}),$ we obtain%
\begin{equation*}
l_{j}^{n}=\dfrac{\pi _{p}}{\lambda _{n}^{2/p}}+\frac{2}{p\lambda _{n}}%
r(z)l_{j}^{n}+\frac{1}{p\lambda _{n}^{2}}q(z)l_{j}^{n}+O\left( \frac{1}{%
\lambda _{n}^{\frac{4+p}{p}}}\right)
\end{equation*}%
or%
\begin{equation*}
q(z)=p\lambda _{n}^{2}\left( \frac{\pi _{p}}{\lambda _{n}^{2/p}l_{j}^{n}}%
\right) \left( \frac{\lambda _{n}^{2/p}l_{j}^{n}}{\pi _{p}}-\frac{2r(z)}{%
p\lambda _{n}}\frac{\lambda _{n}^{2/p}l_{j}^{n}}{\pi _{p}}-1\right) .
\end{equation*}%
Considering (2.6), we can write that for $n\rightarrow \infty $%
\begin{equation*}
\frac{\lambda _{n}^{2/p}l_{j}^{n}}{\pi _{p}}=1,
\end{equation*}%
Then,%
\begin{equation*}
q(x)=\lim_{n\rightarrow \infty }p\lambda _{n}^{2}\left( \frac{\lambda
_{n}^{2/p}l_{j}^{n}}{\pi _{p}}-\frac{2r}{p\lambda _{n}}-1\right) .
\end{equation*}

This completes the proof.

\textbf{Conclusion 2.4.} In the Theorem 2.1., Theorem2.2., Theorem 2.3 and
Theorem 3.1, taking $r(x)=0$ we obtain results of Sturm-Liouville problem
given in \cite{lav}.

\textbf{Conclusion 2.5.} In the Theorem 2.1., Theorem2.2., Theorem 2.3, and
Theorem 3.1, taking $p=2$, we obtain results of inverse nodal problem for
differential pencil.

\end{document}